\documentclass[12pt,leqno,fleqn]{article}
\usepackage{epsf,amsmath,amssymb,amscd,hyperref,srcltx}
\usepackage{amsfonts}
\usepackage{psfrag}
\newcommand{\dps}{\displaystyle}
\newcommand{\id}{\text{\rm id}}
\newcommand{\GL}{\text{\rm GL}}
\newcommand{\T}{^{\sf T}}
\newcommand{\dy}[2]{%
\refstepcounter{equation}%
\LABEL{#1}%
\begin{list}{}{
\topsep 3mm
\leftmargin 18mm
\rightmargin 0cm
\itemsep 0mm
\listparindent 0mm
\parsep 0mm
\itemsep 0mm
\labelsep 0mm
\labelwidth 18mm
}%
\item[\rm (\theequation)\hfill]
#2
\end{list}%
}
\newcommand{\dyz}[1]{%
\refstepcounter{equation}%
\begin{list}{}{
\topsep 3mm
\leftmargin 18mm
\rightmargin 0cm
\itemsep 0mm
\listparindent 0mm
\parsep 0mm
\itemsep 0mm
\labelsep 0mm
\labelwidth 18mm
}%
\item[\rm (\theequation)\hfill]
#1
\end{list}%
}
\newcommand{\dyyz}[1]{\dyz{\raggedright$\dps#1$}}
\newcommand{\dyy}[2]{\dy{#1}{\raggedright$\dps#2$}}
\newcounter{bewering}
\newcommand{\prop}[2]{\refstepcounter{bewering}\vspace{4mm}\noindent{\bf Proposition \thebewering.}\label{#1}{\it #2}}
\newcommand{\propz}[1]{\refstepcounter{bewering}\vspace{4mm}\noindent{\bf Proposition \thebewering.}{\it #1}}
\newcounter{claim}
\newcommand{\clnn}[1]{\refstepcounter{claim}\vspace{4mm}\noindent{\em Claim.}  {\it #1}}
\newcounter{sectie}
\newcommand{\sect}[2]{\refstepcounter{sectie}
\section*{\boldmath \thesectie. #2}%
\label{#1}}
\newcommand{\sectz}[1]{\refstepcounter{sectie}
\section*{\boldmath \thesectie. #1}%
}
\newcommand{\pf}{\vspace{3mm}\noindent{\bf Proof.}\ }
\newcommand{\opf}{\vspace{3mm}\noindent{\em Proof.}\ }
\newcommand{\bx}{\hspace*{\fill} \hbox{\hskip 1pt \vrule width 4pt height 8pt depth 1.5pt \hskip 1pt}

\addvspace{4mm}}
\newcommand{\obx}{\hspace*{\fill} \hbox{$\Box$}

\addvspace{4mm}}
\newcommand{\rf}[1]{{\rm (\ref{#1})}}
\newcommand{\FF}{{\cal F}}
\newcommand{\GG}{{\cal G}}
\newcommand{\MM}{{\cal M}}
\newcommand{\OO}{{\cal O}}
\newcommand{\VV}{{\cal V}}
\newcommand{\LABEL}[1]{\label{#1}}
\newcommand{\rank}{\text{\rm rank}}
\newcommand{\sgn}{\text{\rm sgn}}
\newcommand{\oC}{{\mathbb{C}}}
\newcommand{\oR}{{\mathbb{R}}}
\newcommand{\oZ}{{\mathbb{Z}}}
\renewcommand{\phi}{\varphi}
\newcommand{\Ker}{{\text{\rm Ker}}}
\newcommand{\AS}{\text{\rm AS}}
\newcommand{\IHX}{\text{\rm IHX}}
\newcommand{\gl}{\frak{gl}}
\renewcommand{\phi}{\varphi}
\newcommand{\Space}{((\oC^n)^{\otimes 3})^{C_3}}
\newcommand{\bgcirc}{\raisebox{1.4pt}{\mbox{\scriptsize $\mathbf{\bigcirc}$}}}
\newcommand{\Thetagraaf}{\raisebox{-.25\height}{\scalebox{0.02}{\includegraphics{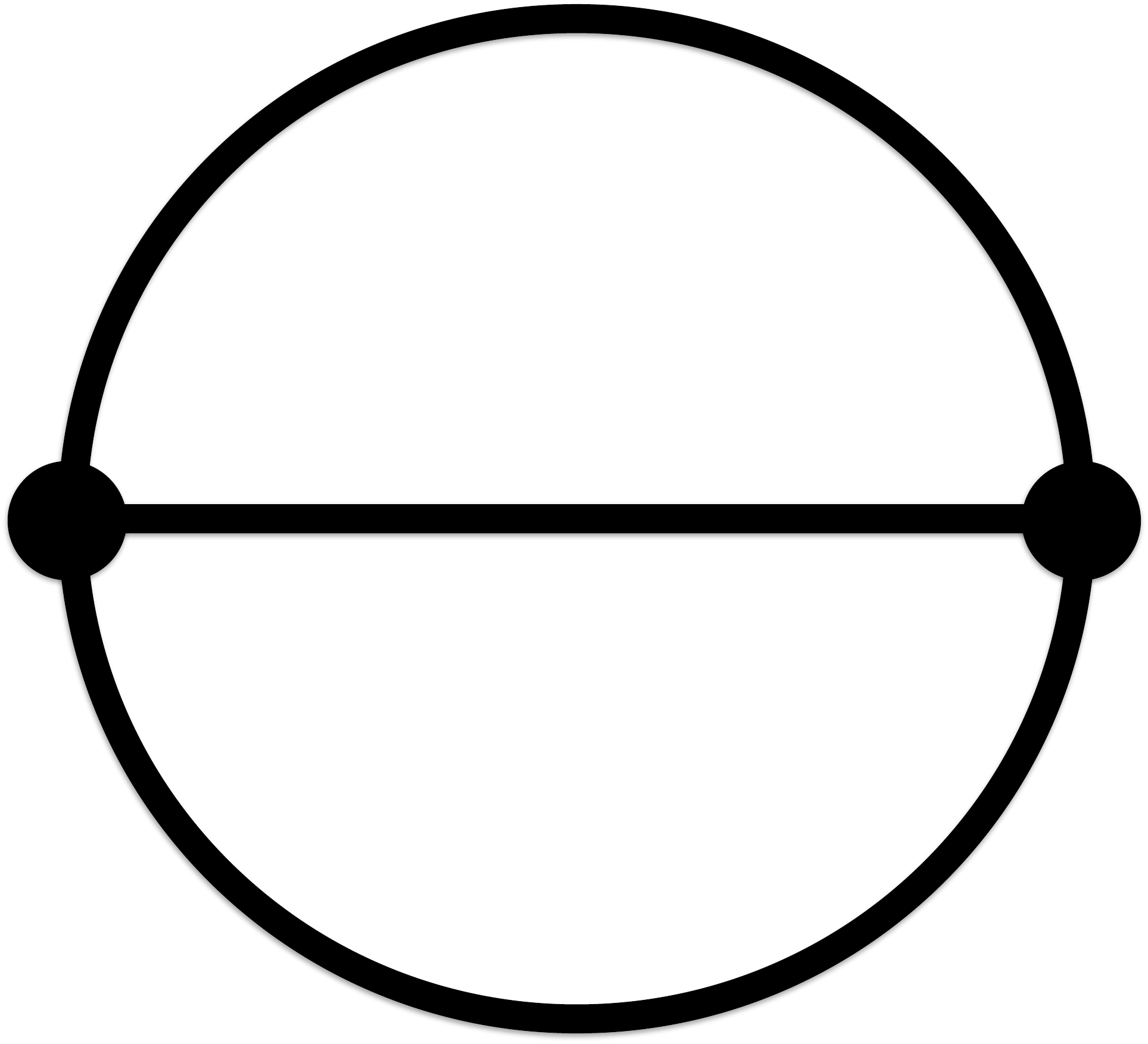}}}}
\begin{document}

\begin{center}
{\large\bf ON LIE ALGEBRA WEIGHT SYSTEMS FOR 3-GRAPHS

}
\vspace{4mm}
Alexander Schrijver \footnote{ University of Amsterdam and CWI, Amsterdam.
Mailing address: CWI, Science Park 123, 1098 XG Amsterdam, The Netherlands.
Email: lex@cwi.nl.
The research leading to these results has received funding from the European Research Council
under the European Union's Seventh Framework Programme (FP7/2007-2013) / ERC grant agreement
n$\mbox{}^{\circ}$ 339109.}

\end{center}

\noindent
{\small{\bf Abstract.}
A {\em $3$-graph} is a connected cubic graph such that each vertex is
is equipped with a cyclic order of the edges incident with it.
A {\em weight system} is a function $f$ on the collection of
$3$-graphs which is {\em antisymmetric}\/:
$f(H)=-f(G)$ if $H$ arises from $G$ by reversing the orientation at one of its vertices, and
satisfies the IHX-equation.
Key instances of weight systems are the functions $\varphi_{\frak{g}}$ obtained from a metric Lie algebra
$\frak{g}$ by taking the structure tensor $c$ of $\frak{g}$ with respect to some
orthonormal basis, decorating each vertex of the $3$-graph by $c$, and contracting along the edges.

We give equations on values of any complex-valued weight system that characterize it as complex Lie algebra
weight system.
It also follows that if $f=\varphi_{\frak{g}}$ for some complex metric Lie algebra $\frak{g}$,
then $f=\varphi_{\frak{g}'}$ for some unique complex reductive metric Lie algebra $\frak{g}'$.
Basic tool throughout is geometric invariant theory.



}

\sectz{Introduction}

A {\em $3$-graph} is a connected, nonempty, cubic graph such that each vertex $v$ is
equipped with a cyclic order of the edges incident with $v$.
Loops and multiple edges are allowed.
Also the `vertexless loop' $\bgcirc$ counts as 3-graph.

$3$-graphs come up in various branches of mathematics, under several names.
They play an important role in studying the Vassiliev knot invariants (see
[6]), and in this context the term 3-graph was introduced by Duzhin, Kaishev, and Chmutov [8].
We adopt this name as it is short and settled in [6].
3-graphs also emerge in the related Chern-Simons topological field theory
(Bar-Natan [2], Axelrod and Singer [1]).
They are in one-to-one correspondence with cubic graphs that are cellularly embedded on a
compact oriented surface, and hence, through graph duality, with triangulations of a
compact oriented surface.
Moreover, $3$-graphs produce a generating set for the algebra of
$O(n,\oC)$-invariant regular functions on the space $\Lambda^3\oC^n$ of
alternating $3$-tensors, by Weyl's `first fundamental theorem' of invariant theory [17].

For the Vassiliev knot invariants, `weight systems' are pivotal.
For 3-graphs they are defined as follows.
Let $\GG$ denote the collection of all 3-graphs, and call
a function $f:\GG\to\oC$ a {\em weight system} if $f$ satisfies the {\em AS-equation}
(for antisymmetry):
$f(H)=-f(G)$ whenever $H$ arises from $G$ by turning the cyclic order of the edges at one of the
vertices of $G$,
and the {\em IHX-equation}\/:
\dy{25se14a}{
$f(~\raisebox{-.2\height}{\scalebox{0.15}{\includegraphics{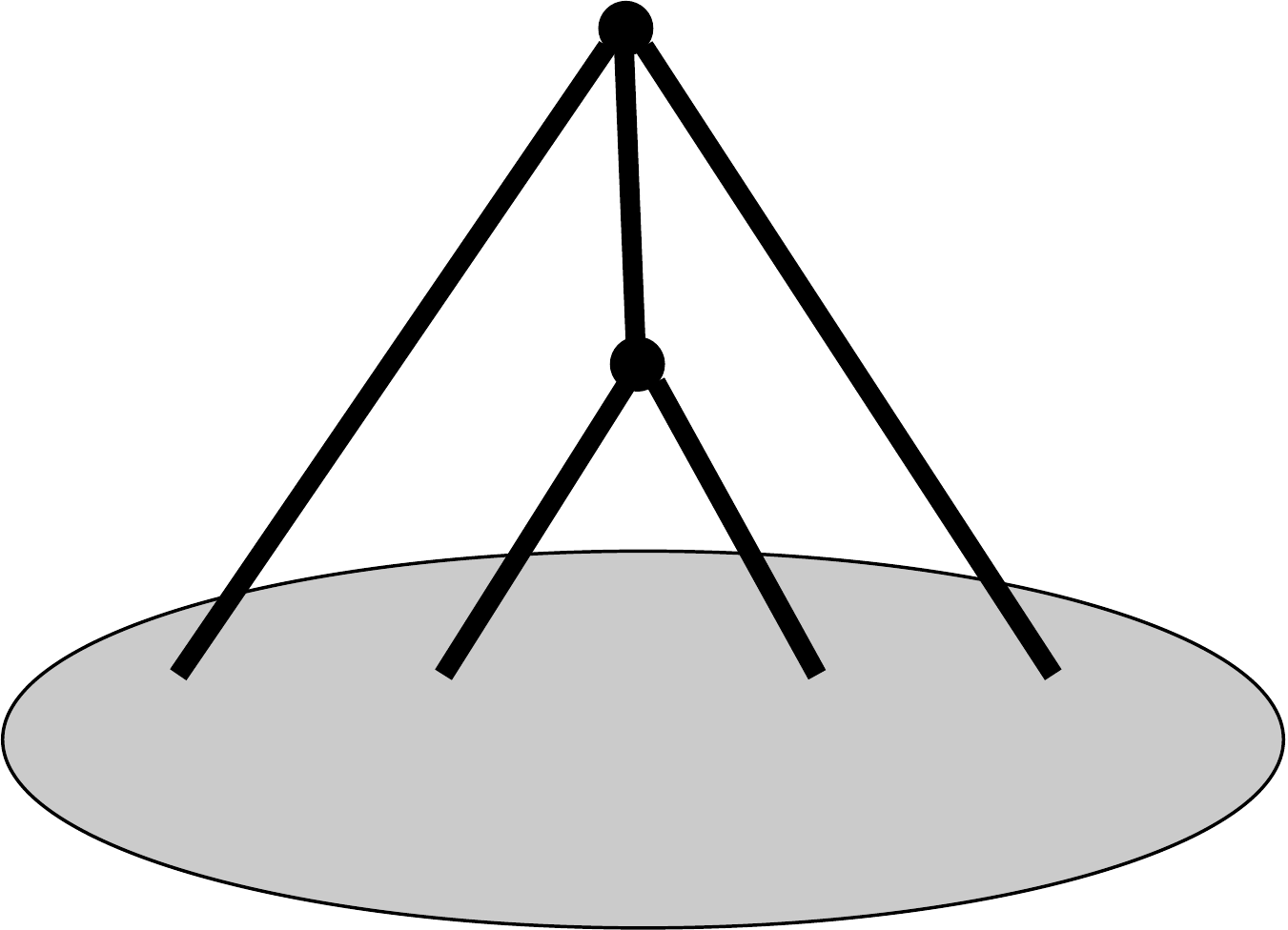}}}~)
~~=~~
f(~\raisebox{-.3\height}{\scalebox{0.15}{\includegraphics{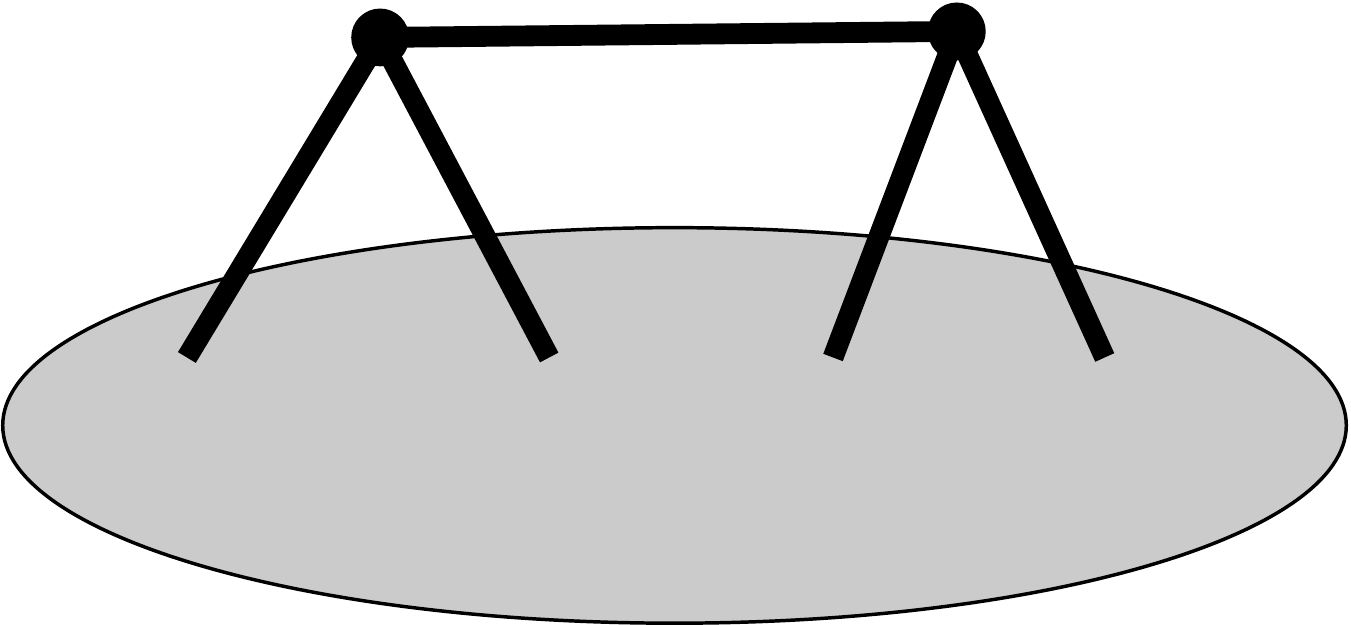}}}~)
~~-~~
f(~\raisebox{-.3\height}{\scalebox{0.15}{\includegraphics{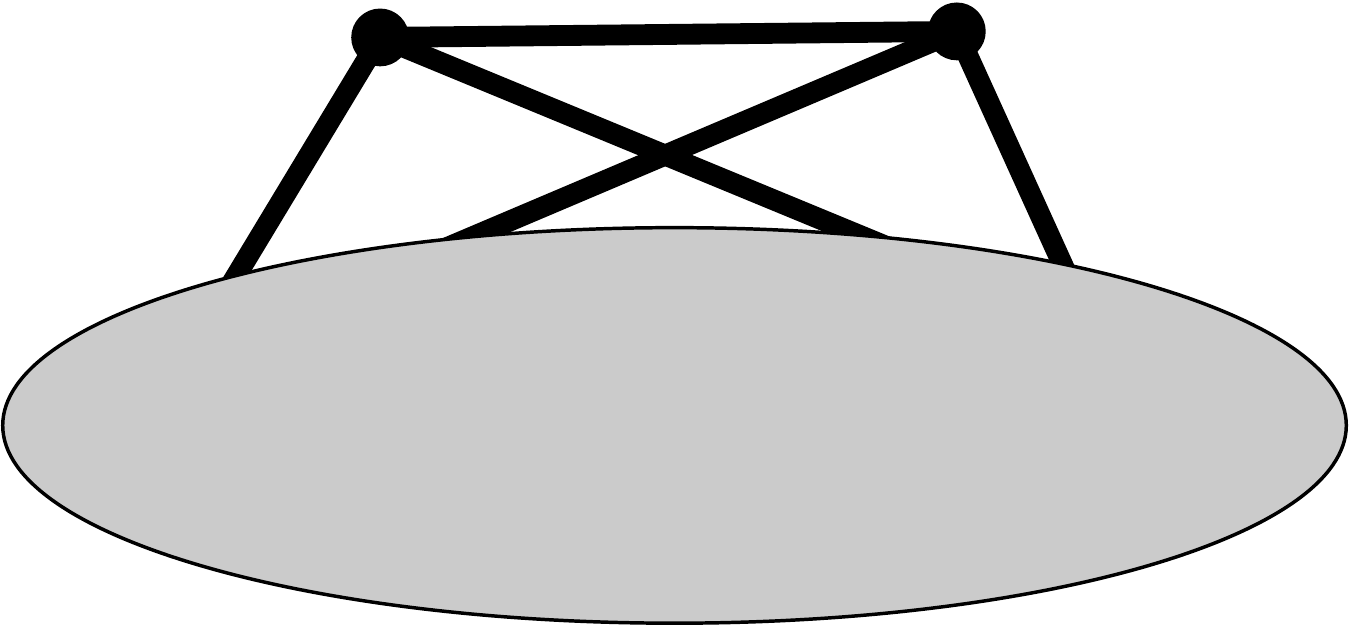}}}~)$.
}
Here the cyclic order of edges at any vertex is given by the clockwise order of edges at the vertex.
The grey areas in \rf{25se14a} represent the remainder of the 3-graphs, the same in each of these areas.

Any $C_3$-invariant tensor $c=(c_{i,j,k})_{i,j,k=1}^n\in((\oC^n)^{\otimes 3})^{C_3}$ gives
the {\em partition function} $p_c:\GG\to\oC$, defined by
\dyy{3ok14a}{
p_c(G):=\sum_{\psi:E(G)\to[n]}\prod_{v\in V(G)}c_{\psi(e_1),\psi(e_2),\psi(e_3)}
}
for any 3-graph $G$, where, for any $v\in V(G)$,
$e_1,e_2,e_3$ denote the edges incident with $v$, in cyclic order.
(As usual, $[k]:=\{1,\ldots,k\}$ for any $k\in\oZ_+=\{0,1,2\ldots\}$.
Moreover, if a group $\Gamma$ acts on a set $X$, then $X^{\Gamma}$ denotes the set of $\Gamma$-invariant
elements of $X$.)

Note that \rf{3ok14a} is invariant under orthonormal transformations of $c$, and that $p_c(\bgcirc)=n$.
($p_c(G)$ is the `partition function' of the `vertex model' $c$, in the sense
of de la Harpe and Jones [11].
It may also be viewed as `edge coloring model' as in Szegedy [15].)

An important class of weight systems is obtained as follows from the structure constants of
metric Lie algebras,
which roots in papers of Penrose [14] and Murphy [13], and the relevance for knot theory was pioneered by
Bar-Natan [2,\linebreak[0]3] and Kontsevich [12].
(In this paper, all Lie algebras are finite-dimensional and complex.)
A {\em metric Lie algebra} is a Lie algebra $\frak{g}$ enriched
with a nondegenerate symmetric bilinear form $\langle.,.\rangle$ which is
{\em ad-invariant}, that is, $\langle[x,y],z\rangle=\langle x,[y,z]\rangle$
for all $x,y,z\in \frak{g}$.
If $b_1,\ldots,b_n$ is an orthonormal basis of $\frak{g}$, the
structure constants $(c_{\frak{g}})_{i,j,k}$ of $\frak{g}$ are characterized by
$(c_{\frak{g}})_{i,j,k}:=\langle[b_i,b_j],b_k\rangle$
for $i,j,k=1,\ldots,n$.
Then $c_{\frak{g}}\in(\frak{g}^{\otimes 3})^{C_3}$,
and $\varphi_{\frak{g}}:=p_{c_{\frak{g}}}$ is a weight system.
Indeed, the set of such structure constants $c_{\frak{g}}$ of metric Lie algebras $\frak{g}$
in $n$ dimensions is equal to the affine variety $\VV_n$ in $\Lambda^3\oC^n$ determined by
the quadratic equations
\dy{12ok14c}{
$\dps\sum_{a=1}^n(x_{i,j,a}x_{a,k,l}+x_{k,i,a}x_{a,j,l}+x_{j,k,a}x_{a,i,l})=0$
\hfill for $i,j,k,l=1,\ldots,n$.
}
This directly implies that $p_x$ satisfies the AS- and IHX-equations;
that is, $p_x$ is a weight system.

As was shown by Bar-Natan [4],
the functions $\varphi_{\frak{sl}(n)}$ connect to basic graph theory properties like
edge-colorability and planarity, and the four-color theorem can be expressed as a relation between the
zeros and the degree of $\varphi_{\frak{sl}(n)}(G)$ (as polynomial in $n$). 

It is easy to construct a weight system that is not equal to $\varphi_{\frak{g}}$ for any
metric Lie algebra: just take a different Lie algebra for each number of vertices of the 3-graph $G$.
More interestingly, Vogel [16] showed that there is a weight system that is no Lie algebra weight system
even when restricted to 3-graphs with 34 vertices.

The AS and IHX conditions are `2-term' and `3-term' relations.
We give a characterization of those weight systems that come from an $n$-dimensional Lie algebra by adding
an `$(n+1)!$-term' relation.
To illustrate it, for $n=2$ it is the following 6-term relation:
\dyz{
$
f(~\raisebox{-.2\height}{\scalebox{0.08}{\includegraphics{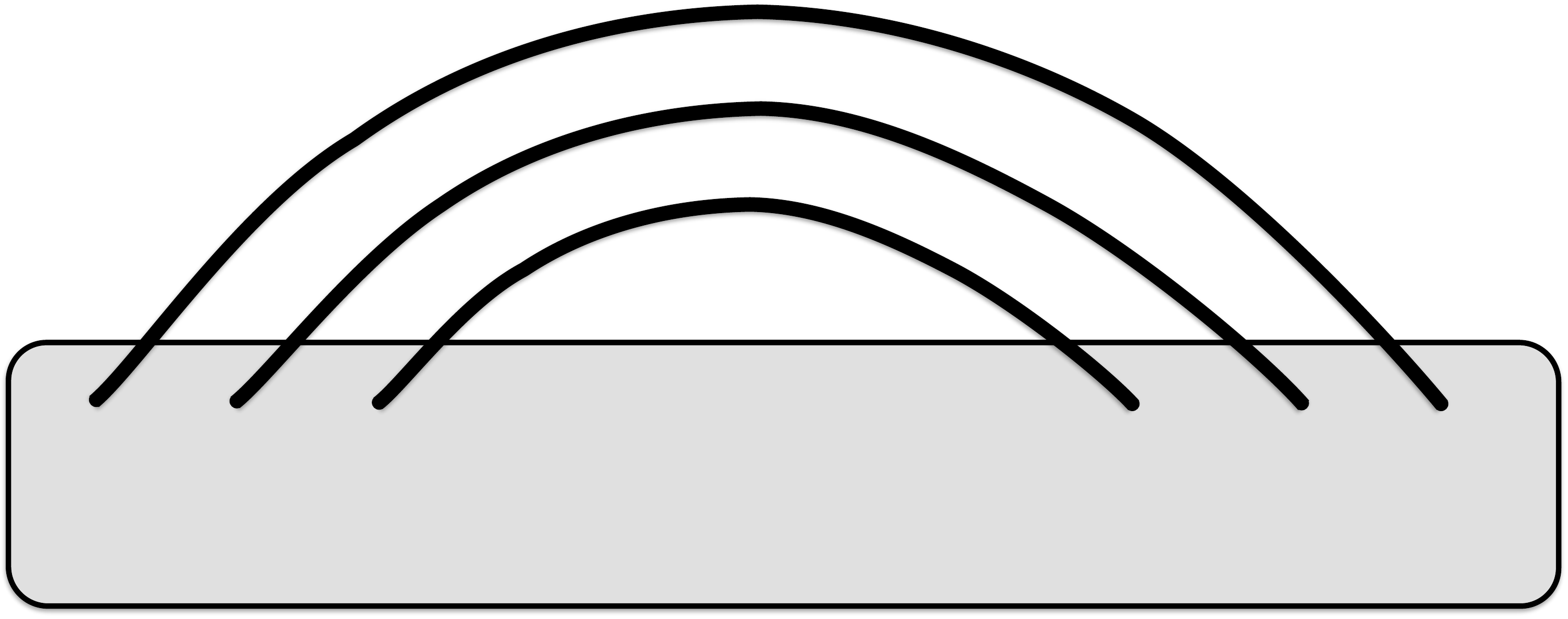}}}~)
~~+~~
f(~\raisebox{-.2\height}{\scalebox{0.08}{\includegraphics{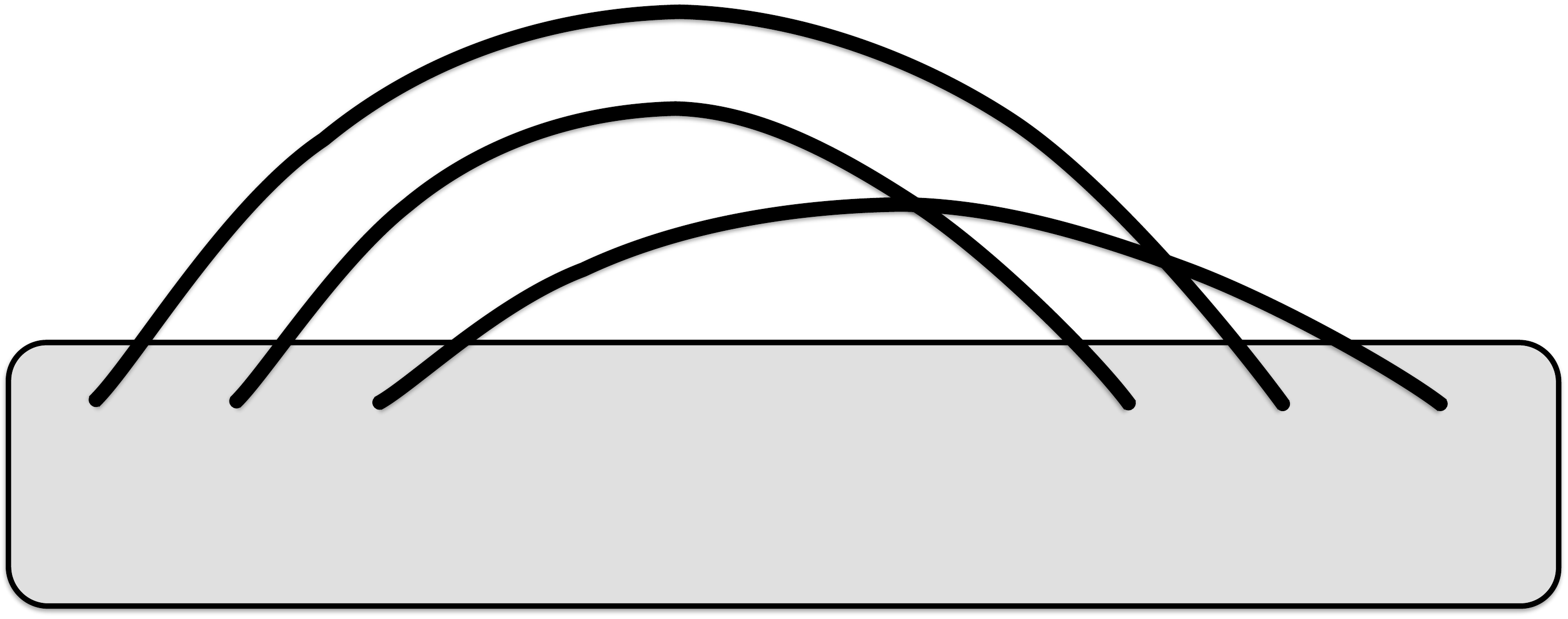}}}~)
~~+~~
f(~\raisebox{-.2\height}{\scalebox{0.08}{\includegraphics{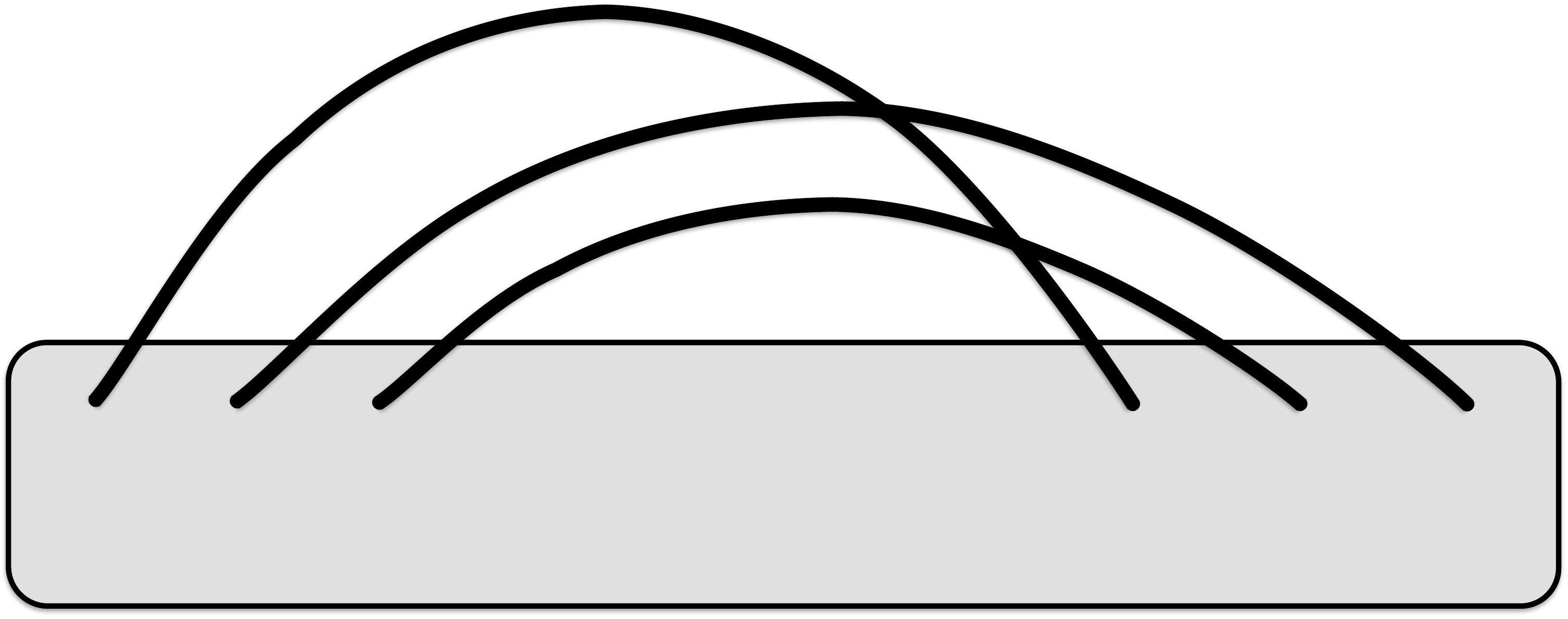}}}~)
$\\
$=~~
f(~\raisebox{-.2\height}{\scalebox{0.08}{\includegraphics{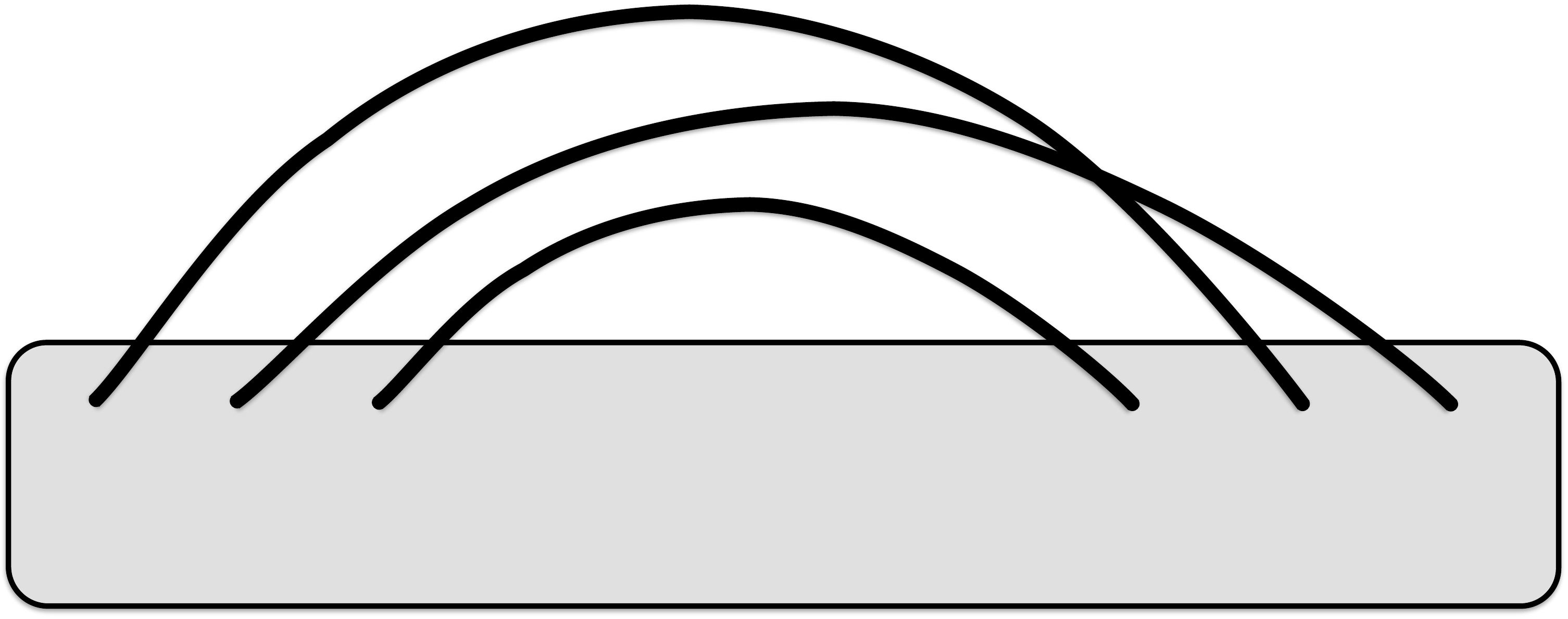}}}~)
~~+~~
f(~\raisebox{-.2\height}{\scalebox{0.08}{\includegraphics{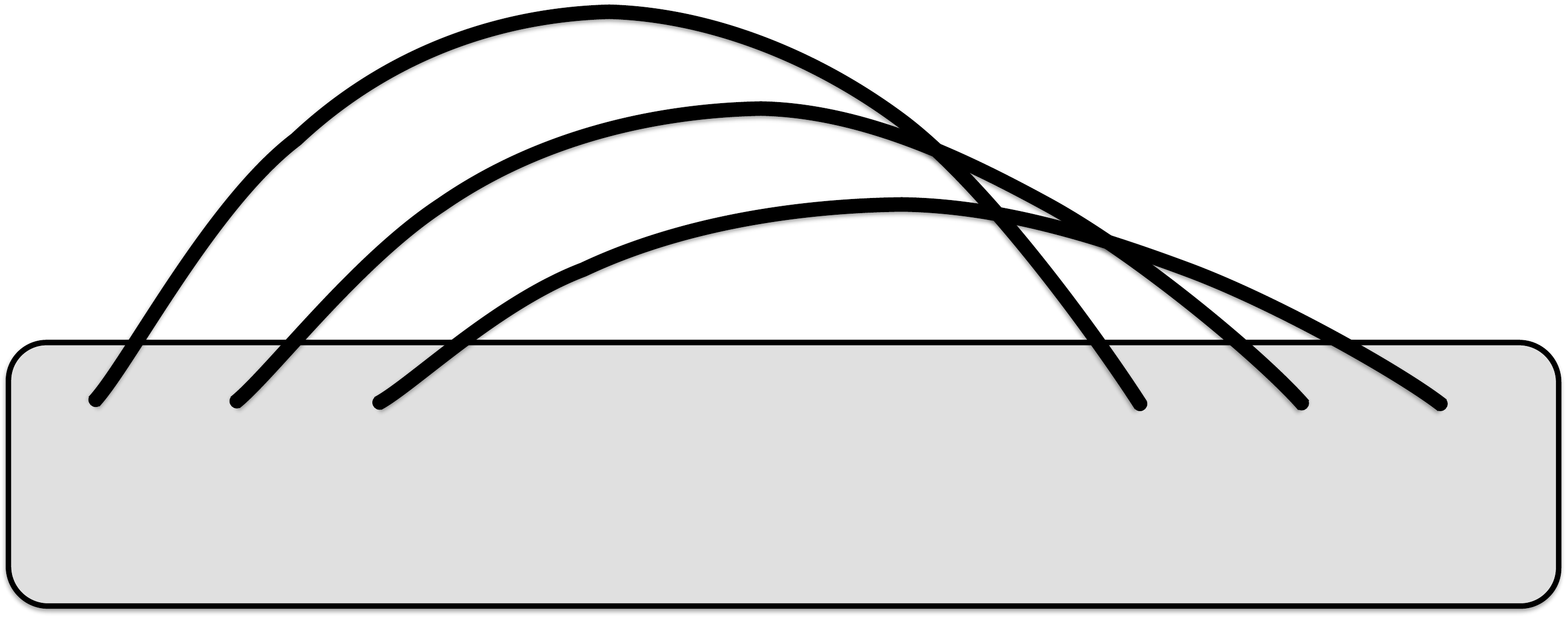}}}~)
~~+~~
f(~\raisebox{-.2\height}{\scalebox{0.08}{\includegraphics{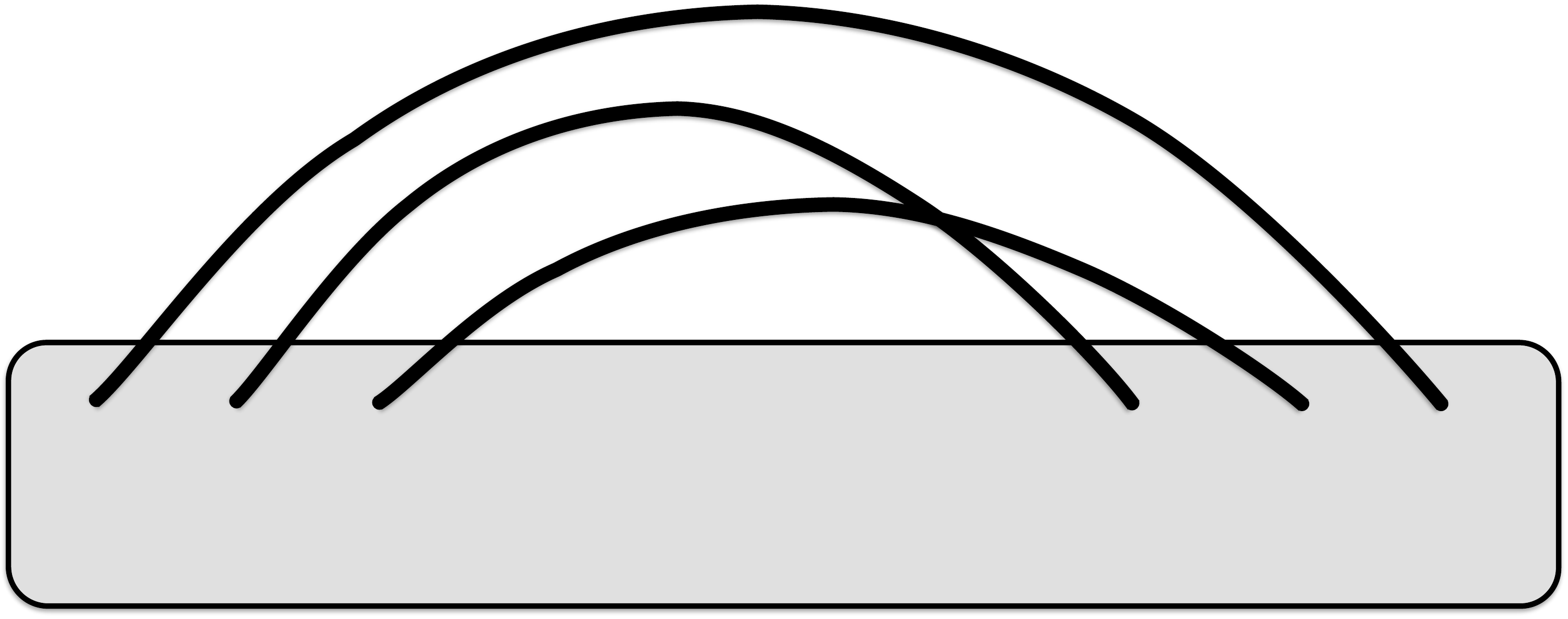}}}~)
$~.
}
To describe it in general, we need the following notion (cf.\ [6]).
A {\em $k$-legged fixed diagram} is a graph with trivalent vertices,
each equipped with a cyclic order of the edges incident with it, and moreover precisely $k$
univalent vertices (called {\em legs}), labeled $1,\ldots,k$.
(Connectivity is not required.)
Let $\FF_k$ denote the collection of all $k$-legged fixed diagrams.

So $\FF_0$ is the collection of disjoint unions of 3-graphs.
Any function $f$ on 3-graphs can be extended to $\FF_0$ by the `multiplicativity rule'
$f(G\sqcup H)=f(G)f(H)$, where
$\sqcup$ denotes disjoint union, and setting $f(\emptyset):=1$.

For $G,H\in\FF_k$, let $G\cdot H$ be the graph in $\FF_0$ obtained from the disjoint union
of $G$ and $H$ by
identifying the $i$-labeled legs in $G$ and $H$ and joining the two incident edges to one edge
(thus forgetting these end vertices as vertex), for each $i=1,\ldots,k$.

For $\pi\in S_k$, let $P_{\pi}$ be the $2k$-legged fixed diagram consisting of $k$ disjoint
edges $e_1,\ldots,e_k$, where the ends of edge $e_i$ are labeled $i$ and $k+\pi(i)$ (for
$i=1,\ldots,k$).
Then we call the following $k!$-term relation the {\em $\Delta_k$-equation}:
\dy{21ok14b}{
$\dps\sum_{\pi\in S_k}\sgn(\pi)f(P_{\pi}\cdot H)=0$\hfill for each $2k$-legged fixed diagram $H$.
}

\noindent
{\bf Theorem.}
{\it
Let $f:\GG\to\oC$ be a weight system.
Then there exists a complex reductive metric Lie algebra $\frak{g}$ with $f=\varphi_{\frak{g}}$
if and only if $f(\bgcirc)\in\oZ_+$ and $f$ satisfies the $\Delta_{f(\bigcirc)+1}$-equation.
If $\frak{g}$ exists, it is unique.
}

\bigskip
Although \rf{21ok14b} may look like a linear constraint separately for each fixed number of vertices of 3-graphs
(as the AS- and IHX-equations are), it
in fact interconnects 3-graphs with different numbers of vertices, since $P_{\pi}\cdot H$ can be a
disjoint union of 3-graphs, taking multiplicativity of $f$ as above.
So \rf{21ok14b} is a polynomial relation between $f$-values of 3-graphs.
It describes $f$ as common zero of a set of polynomials in the ring $\oC[\GG]$ formally
generated by the collection $\GG$ of 3-graphs (which are connected by definition), in which the
disjoint union $\sqcup$ is taken as multiplication.
We note that, for any $k\in\oZ_+$,
the $\Delta_k$-equation implies that $f(\bgcirc)$ is a nonnegative integer strictly less than $k$
(by taking $H:=P_{\id}$, where $\id$ is the identity permutation in $S_k$).

The theorem implies that if $f=\phi_{\frak{g}}$ for some metric Lie algebra $\frak{g}$,
then $f=\phi_{\frak{g}}$ for a unique reductive metric Lie algebra $\frak{g}$.
Indeed, let $f=\phi_{\frak{g}}$ for some $n$-dimensional Lie algebra $\frak{g}$.
So $f=p_{c_{\frak{g}}}$ and $f(\bgcirc)=n$.
We show that $f$ satisfies the $\Delta_{n+1}$-equation.
Take a $2(n+1)$-legged fixed diagram $H$ and consider formulas \rf{3ok14a} and \rf{21ok14b} for
$f:=p_{c_{\frak{g}}}$.
The summations over $\pi$ and $\psi$ can be interchanged.
For each fixed $\psi:E(H)\to[n]$, we need to add up $\sgn(\pi)$ over all those $\pi\in S_{n+1}$
for which, for each $i\in[n+1]$, legs $i$ and $k+\pi(i)$ of $H$ have the same $\psi$-value.
As $n<n+1$, there exist two legs $i,j\in[n+1]$ of $H$ with the same $\psi$-value.
Let $\sigma\in S_{n+1}$ be the transposition of $i$ and $j$.
Then we can pair up each $\pi$ with $\pi\sigma$, and in the summation they cancel.
So for each $\psi$, the sum is 0, and therefore the $\Delta_{n+1}$-equation holds.

These arguments also yield the necessity of the condition in the theorem.
We prove sufficiency in Section \ref{8ok14l} and uniqueness in Section \ref{8ok14j}.

By the theorem, if the $\Delta_{f(\bigcirc)+1}$-equation holds, then
there exist unique (up to permuting indices) simple Lie algebras $\frak{g}_1,\ldots\frak{g}_t$ and
nonzero complex numbers $\lambda_1,\ldots,\lambda_t$ such that
\vspace*{-4pt}
\dyyz{
f(G)=\sum_{i=1}^t\lambda_i^{\frac12|V(G)|}\varphi_{\frak{g}_i}(G)\vspace*{-4pt}
}
for each $3$-graph $G\neq\bgcirc$, taking the Killing forms as metrics.
So any linear combination of 3-graphs `detected' (to be nonzero) by a Lie algebra weight system,
is detected by a simple Lie algebra weight system.

It can be proved that if $n:=f(\bgcirc)\in\oZ_+$, then the $\Delta_{n+1}$-equation can be
replaced by the equivalent condition that for each $k\in\oZ_+$, the rank of the
$\FF_k\times\FF_k$ matrix $C_{f,k}:=(f(G\cdot H)_{G,H\in\FF_k})$ is at most $n^k$.
A weaker, but also equivalent condition is that there exists an $m\in\oZ_+$ such that the rank of
$C_{f,2(n+1)m}$ is less than the dimension of the space of all $\GL(d)$-invariant tensors in
$\gl(d)^{\otimes m}$, where $d:=n^{n+1}+1$.

Our proof is based on some basic theorems from invariant theory
(Weyl's first and second fundamental
theorem for the orthogonal group, and the unique closed orbit theorem;
cf.\ [5],[10]), and roots in methods used in [7], [9], and
[15].
For any $n\in\oZ_+$ and any 3-graph $G$, let $p(G)$ be the regular function on $((\oC^n)^{\otimes 3})^{C_3}$ defined by
\dy{20ok14a}{
$p(G)(x):=p_x(G)$~~~~~for $x\in((\oC^n)^{\otimes 3})^{C_3}$
}
(cf.\ \rf{3ok14a}).
Then each $p(G)$ is $O(n,\oC)$-invariant, and the
first fundamental theorem of invariant theory
implies that the algebra of $O(n,\oC)$-invariant
regular functions on $((\oC^n)^{\otimes 3})^{C_3}$ is generated by
$\{p(G)\mid G$ 3-graph$\}$.

Note that $O(n,\oC)$ acts naturally on the affine variety $\VV_n$ defined by \rf{12ok14c}, and that
for any two metrized Lie algebras $\frak{g}$ and $\frak{g}'$ one trivially has:
$\frak{g}=\frak{g}'$ if and only if
$c_{\frak{g}}$ and $c_{\frak{g}'}$ belong to the same $O(n,\oC)$-orbit on $\VV_n$.
Moreover, the closed orbit theorem implies that
$\varphi_{\frak{g}}=\varphi_{\frak{g}'}$ if and only if the closures of the orbits 
$O(n,\oC)\cdot c_{\frak{g}}$ and $O(n,\oC)\cdot c_{\frak{g}'}$ intersect; that is,
if and only if they project to the same point in $\VV_n//O(n,\oC)$.

The proof implies that a metric Lie algebra $\frak{g}$ is reductive
if and only if the orbit $O(\frak{g})\cdot c_{\frak{g}}$ is closed.
Hence, for each $n$, there is a one-to-one correspondence between the points
in the orbit space $\VV_n//O(n,\oC)$ and the $n$-dimensional complex reductive metric Lie algebras.

\sect{8ok14l}{Proof of the theorem: existence of $\frak{g}$}

Let $f:\GG\to\oC$ satisfy the $\Delta_{n+1}$-equation \rf{21ok14b},
where $n:=f(\bgcirc)\in\oZ_+$.
As above, we extend $f$ to the collection $\FF_0$ of disjoint unions of 3-graphs by the rule that
$f(\emptyset)=1$ and $f(G\sqcup H)=f(G)f(H)$ for all $G,H\in\FF_0$
(where $\sqcup$ denotes disjoint union).
For any $k$,
let $\oC\FF_k$ be the linear space of formal $\oC$-linear combinations of elements of $\FF_k$.
Any (bi-)linear function on $\FF_k$ can be extended (bi-)linearly to $\oC\FF_k$.
Taking $\sqcup$ as product, $\oC\FF_0$ becomes an algebra (which is equal to
$\oC[\GG]$ described above), and $f$ becomes an algebra homomorphism
$\oC\FF_0\to\oC$.
Similarly, $p$ (as defined in \rf{20ok14a}) extends to an algebra homomorphism $\oC\FF_0\to\OO(\Space)$.
(As usual, $\OO(\cdot)$ denotes the algebra of $\oC$-valued regular functions on $\cdot$.)

\prop{12me14d}{
$\Ker(p)\subseteq\Ker(f)$.
}

\pf
Let $\gamma\in\oC\FF_0$ with $p(\gamma)=0$.
We prove that $f(\gamma)=0$.
As each homogeneous component of $p(\gamma)$ is 0,
we can assume that $\gamma$ is a linear combination of graphs in $\FF_0$ that all have the same number
of vertices, say $k$ (which is necessarily even).

Let $\MM$ be the collection of perfect matchings on $[3k]$.
We can naturally identify $\MM$ with
the set of $3k$-legged fixed diagrams with no trivalent vertices and no copies of $\bgcirc$.

Let $H$ be the the $3k$-legged fixed diagram with precisely $k$ trivalent vertices $v_1,\ldots,v_k$
and $3k$ legs,
where $v_i$ is adjacent to legs $3i-2$, $3i-1$, $3i$, in order.
(So $H$ is the disjoint union of $k$ copies of the tri-star $K_{1,3}$.)
Then each graph in $\FF_0$ with $k$ vertices is equal to 
$M\cdot H$ for at least one $M\in\MM$.
Hence we can write
\dyyz{
\gamma=\sum_{M\in\MM}\lambda(M)M\cdot H
}
for some $\lambda:\MM\to\oC$.

The symmetric group $S_{3k}$ acts naturally on $\FF_{3k}$ (by permuting leg-labels).
Let $Q$ be the group of permutations $\sigma\in S_{3k}$ with $H^{\sigma}=H$.
(So $Q$ is the wreath product of the cyclic group $C_3$ with $S_k$.
It stabilizes the partition $\{\{3i-2,3i-1,3i\}\mid i=1,\ldots,k\}$ of $[3k]$ and permutes each
class in this partition cyclically.)
Since $M^{\sigma}\cdot H=M\cdot H^{\sigma^{-1}}=M\cdot H$ for each $M\in\MM$ and $\sigma\in Q$,
we can assume that $\lambda$ is invariant under the action of $Q$ on $\MM$.

Define linear functions $F_M$ (for $M\in\MM$) and $F$ on $(\oC^n)^{\otimes 3k}$ by
\dyy{7se14c}{
F_M(a_1\otimes\cdots\otimes a_{3k}):=\prod_{ij\in M}a_i\T a_j
\text{~~~and~~~}
F:=\sum_{M\in\MM}\lambda(M)F_M,
}
for $a_1,\ldots,a_{3k}\in\oC^n$.
($ij$ stands for the unordered pair $\{i,j\}$; so $ij=ji$.)
Note that $F_M(x^{\otimes k})=p(M\cdot H)(x)$ for any $x\in\Space$.
Hence $F(x^{\otimes k})=p(\gamma)(x)=0$.
We show that this implies that $F=0$.

Indeed, suppose $F(u_1\otimes\cdots\otimes u_k)\neq 0$ for some $u_1,\ldots,u_k\in(\oC^n)^{\otimes 3}$.
Since $F$ is $Q$-invariant (as $\lambda$ is $Q$-invariant), we can assume that each $u_i$ is $C_3$-invariant.
For $y\in\oC^k$, define the $C_3$-invariant tensor $b_y:=y_1u_1+\cdots+y_ku_k$.
As $F$ is $Q$-invariant,
the coefficient of the monomial $y_1\cdots y_k$ in the polynomial $F(b_y^{\otimes k})$ is equal
to $k!\cdot F(u_1\otimes\cdots\otimes u_k)\neq 0$.
So the polynomial is nonzero, hence $F(b_y^{\otimes k})\neq 0$ for some $y\in\oC^k$, a contradiction.
Therefore, $F=0$.

Define for each multiset $N$ of singletons and unordered pairs from $[3k]$ the monomial $q_N$ on $S^2\oC^{3k}$
(= the set of symmetric matrices in $\oC^{3k\times 3k}$), 
and define moreover the polynomial $q$ on $S^2\oC^{3k}$ by:
\dyy{7se14b}{
q_N(X):=\prod_{ij\in N}X_{i,j}
\text{~~~and~~~}
q:=\sum_{M\in\MM}\lambda(M)q_M,
}
for $X=(X_{i,j})\in S^2\oC^{3k}$.
Note that for each monomial $\mu$ on $S^2\oC^{3k}$ there is a unique multiset $N$ of singletons and unordered pairs
from $[3k]$ with $\mu=q_N$.
Now $F=0$ implies
\dy{6se14a}{
$q(X)=0$ if $\rank(X)\leq n$.
}
Indeed, if $\rank(X)\leq n$, then there exist $a_1,\ldots,a_{3k}\in\oC^n$ such that
$X_{i,j}=a_i\T a_j$ for all $i,j=1\ldots,3k$.
By \rf{7se14c} and \rf{7se14b},
$q(X)=F(a_1\otimes\cdots\otimes a_{3k})=0$,
proving \rf{6se14a}.

By the second fundamental theorem of invariant theory
(cf.\ [10] Theorem 12.2.12),
\rf{6se14a} implies that $q$ belongs to the ideal in $\OO(S^2\oC^{3k})$ generated by the
$(n+1)\times (n+1)$ minors of $X\in S^2\oC^{3k}$.
That is, $q$ is a linear combination of polynomials $\det(X_{I,J})q_N(X)$, where $I,J\subseteq[3k]$
with $|I|=|J|=n+1$ and where $N$ is a multiset of singletons and unordered pairs from $[3k]$.
Here $X_{I,J}$ denotes the $I\times J$ submatrix of $X$.

Now such triples $I,J,N$ occur in two kinds: (1) those with $I\cap J=\emptyset$ and $N$ a perfect matching
on $[3k]\setminus(I\cup J)$, in which case {\em all} monomials occurring in $\det(X_{I,J})q_N(X)$
are equal to $q_M$ for some $M\in\MM$;
and (2) all other triples $I,J,N$, in which case {\em none} of the monomials occurring in
$\det(X_{I,J})q_N(X)$ is equal to $q_M$ for some $M\in\MM$.
Since $q(X)$ consists completely of monomials $q_M$ with $M\in\MM$, we can ignore all triples
of kind (2), and conclude that $q(X)$ is a linear combination of $\det(X_{I,J})q_N(X)$ with
$I,J,N$ of kind (1).

For any $M\in\MM$, define $\Gamma(q_M):=M\cdot H$, and extend $\Gamma$ linearly to linear combinations of
the $q_M$ for $M\in\MM$.
Then $\Gamma(q)=\gamma$.
Moreover, for each $I,J,N$ of kind (1), by the $\Delta_{n+1}$-equation for $f$,
$f(\Gamma(\det(X_{I,J})q_N(X)))=0$.
As $\gamma$ is a linear combination of elements $\Gamma(\det(X_{I,J})q_N(X))$,
we have $f(\gamma)=0$, as required.
\bx

By this proposition,
there exists a linear function $\Phi:p(\oC\FF_0)\to\oC$ such that $\Phi\circ p=f$.
Then $\Phi$ is an algebra homomorphism, since for $G,H\in\FF_0$ one has
$\Phi(p(G)p(H))=\Phi(p(G\sqcup H))=f(G\sqcup H)=f(G)f(H)=\Phi(p(G))\Phi(p(H))$.

By the first fundamental theorem of invariant theory,
\dyy{15se14b}{
\OO(\Space)^{O(n)}=p(\oC\FF_0)
}
(setting $O(n):=O(n,\oC)$).
So $\Phi$ is an algebra homomorphism $\OO(\Space)^{O(n)}\to\oC$.
Hence the affine $O(n)$-variety
\dy{15se14c}{
$\VV:=\{x\in \Space\mid q(x)=\Phi(q)$ for each $q\in\OO(\Space)^{O(n)}\}$
}
is nonempty (as $O(n)$ is reductive).
By \rf{15se14b} and by substituting $q=p(G)$ in \rf{15se14c},
\dyyz{
\VV:=\{x\in \Space\mid p_x=f\}.
}
Hence as $\VV\neq\emptyset$ there exists $c\in\Space$ with $p_c=f$.
We choose $c$ such that the orbit $O(n)\cdot c$ is a closed.
This is possible by the unique closed orbit theorem (cf.\ Brion [5]), which also implies that
$c$ is contained in each nonempty $O(n)$-invariant closed subset of $\VV$.

Then $c$ gives the required Lie algebra:

\propz{
$c=c_{\frak{g}}$ for some complex reductive metric Lie algebra $\frak{g}$.
}

\pf
We extend $p(G)$ to a function $\widehat p$ on fixed diagrams as follows.
For each $k$ and $G\in\FF_k$, let $\widehat p(G):\Space\to(\oC^n)^{\otimes k}$ be defined
by
\dyyz{
\widehat p(G)(x):=\sum_{\psi:E(G)\to[n]}
\big(
\prod_{v\in V_3(G)}x_{\psi(e_1),\psi(e_2),\psi(e_3)}
\big)
\bigotimes_{j=1}^kb_{\psi(\varepsilon_j)}
}
for $x\in\Space$,
where $V_3(G)$ is the set of trivalent vertices of $G$, $e_1,e_2,e_3$ are the edges incident with
$v$, in order, and $\varepsilon_j$ is the edge incident with leg labeled $j$ (for $j=1,\ldots,k$).
Moreover, $b_1,\ldots,b_n$ is the standard basis of $\oC^n$.

Then for all $G,H\in\FF_k$ and $x\in\Space$,
\dyy{21ok14a}{
\widehat p(G)(x)\cdot \widehat p(H)(x)=p(G\cdot H)(x),
}
where $\cdot$ denotes the standard inner product on $(\oC^n)^{\otimes k}$.

\clnn{
For each $k$ and $\tau\in\oC\FF_k$, if $f(\tau\cdot H)=0$ for each $H\in\FF_k$,
then $\widehat p(\tau)(c)=0$.
}

\opf
As $\widehat p(\tau)$ is $O(n)$-equivariant, 
it suffices to show that $\widehat p(\tau)$ has a zero $x$ in $\VV$, since
then the $O(n)$-stable closed set $\{x\in\VV\mid \widehat p(\tau)(x)=0\}$ is nonempty,
and hence must contain $c$.

Suppose that no such zero exists.
Then the functions $\widehat p(\tau)$ and $p(G)-f(G)$ (for $G\in\GG$)
have no common zero.
Hence, by the Nullstellensatz, there exist regular functions
$s:\Space\to(\oC^n)^{\otimes k}$ and $g_1,\ldots,g_t:\Space\to\oC$, and
$G_1,\ldots,G_t\in\GG$ such that
\dyyz{
\widehat p(\tau)(x)\cdot s(x)+\sum_{i=1}^t(p(G_i)(x)-f(G_i))g_i(x)=1
}
for all $x\in\Space$.
Now $\widehat p(\tau)$ and $p(G_1),\ldots,p(G_t)$ are $O(n)$-equivariant.
Hence, by applying the Reynolds operator, we can assume that
also $s$ and $g_1,\ldots,g_k$ are $O(n)$-equivariant.
Then by the first fundamental theorem of invariant theory, $s=\widehat p(\beta)$ for some $\beta\in\oC\FF_k$,
and $g_i=p(\gamma_i)$ for some $\gamma_i\in\oC\FF_0$, for $i=1,\ldots,t$.
This gives, with \rf{21ok14a},
\dyyz{
1=\widehat p(\tau)(c)\cdot\widehat p(\beta)(c)+\sum_{i=1}^t(p(G_i)(c)-f(G_i))p(\gamma_i)(c)
=
p(\tau\cdot\beta)(c)+\sum_{i=1}^t(p(G_i)(c)-f(G_i))p(\gamma_i)(c)
=f(\tau\cdot\beta)+\sum_{i=1}^t(f(G_i)-f(G_i))f(\gamma_i)=0,
}
a contradiction, proving the Claim.
\obx

Let $\AS\in\oC\FF_3$ and $\IHX\in\oC\FF_4$ be extracted from the the AS-
and IHX-equations; that is,
\dyyz{
\AS:=
~~
\raisebox{-.4\height}{\scalebox{0.1}{\includegraphics{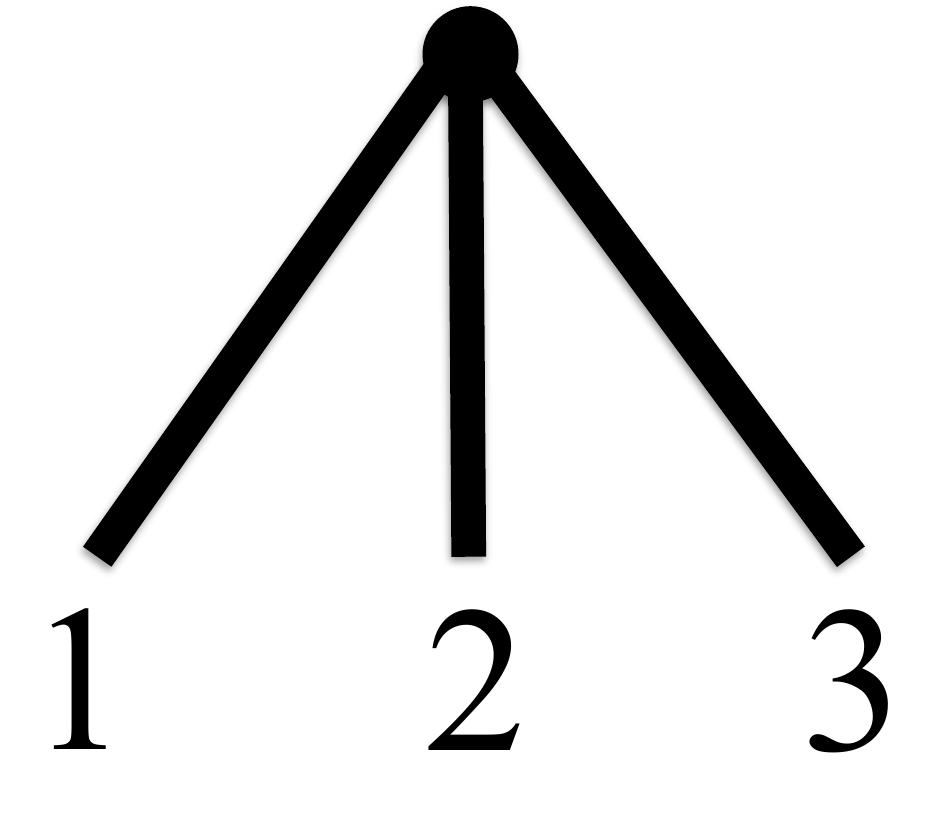}}}~~
+
~~
\raisebox{-.4\height}{\scalebox{0.1}{\includegraphics{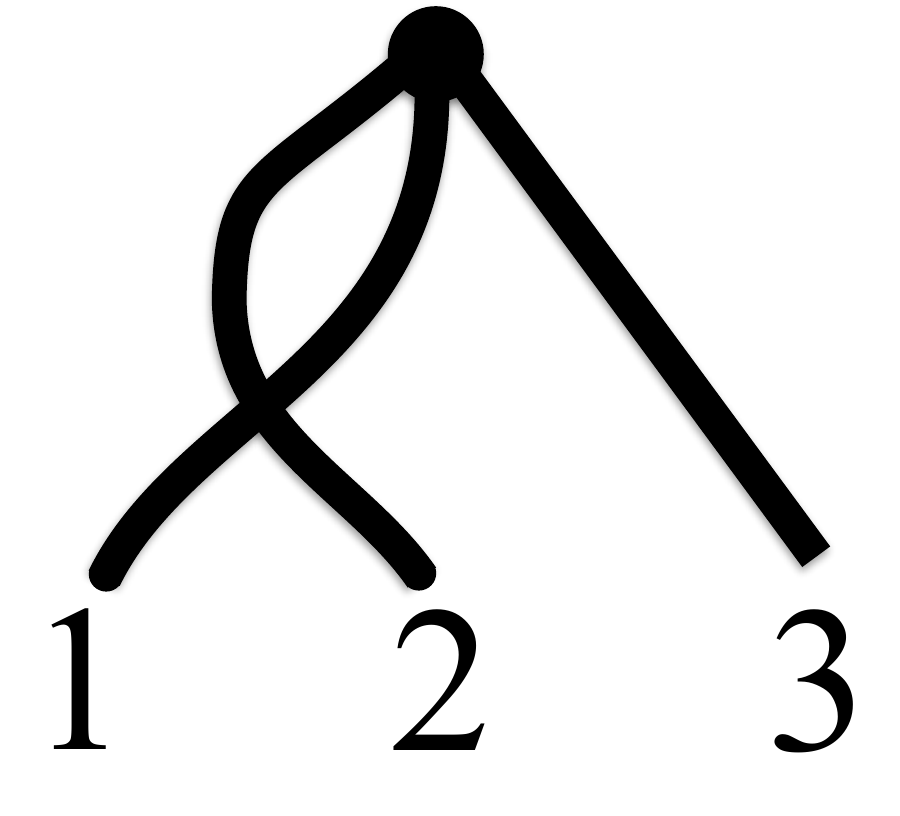}}}~~,
~~~
\IHX:=
~~
\raisebox{-.4\height}{\scalebox{0.1}{\includegraphics{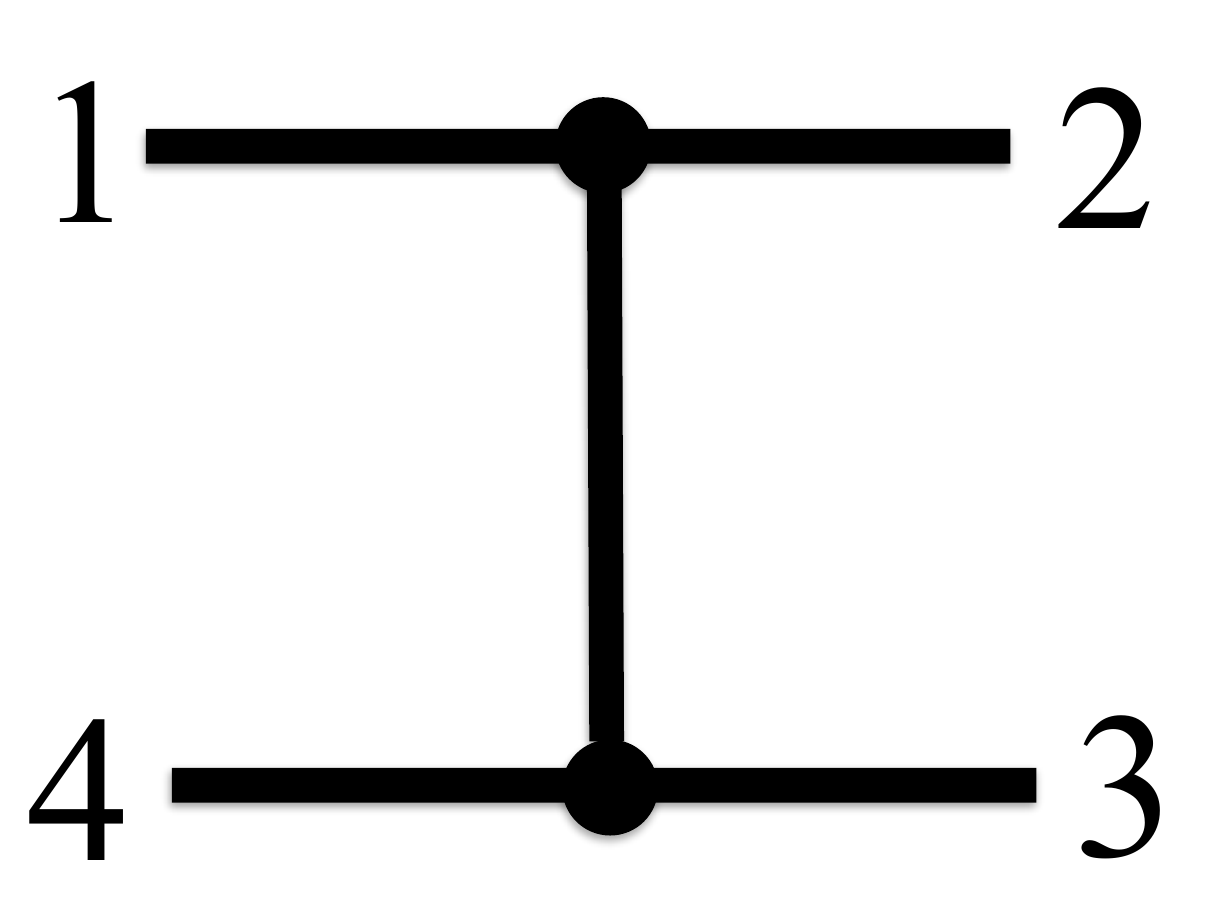}}}~~
-
~~
\raisebox{-.4\height}{\scalebox{0.1}{\includegraphics{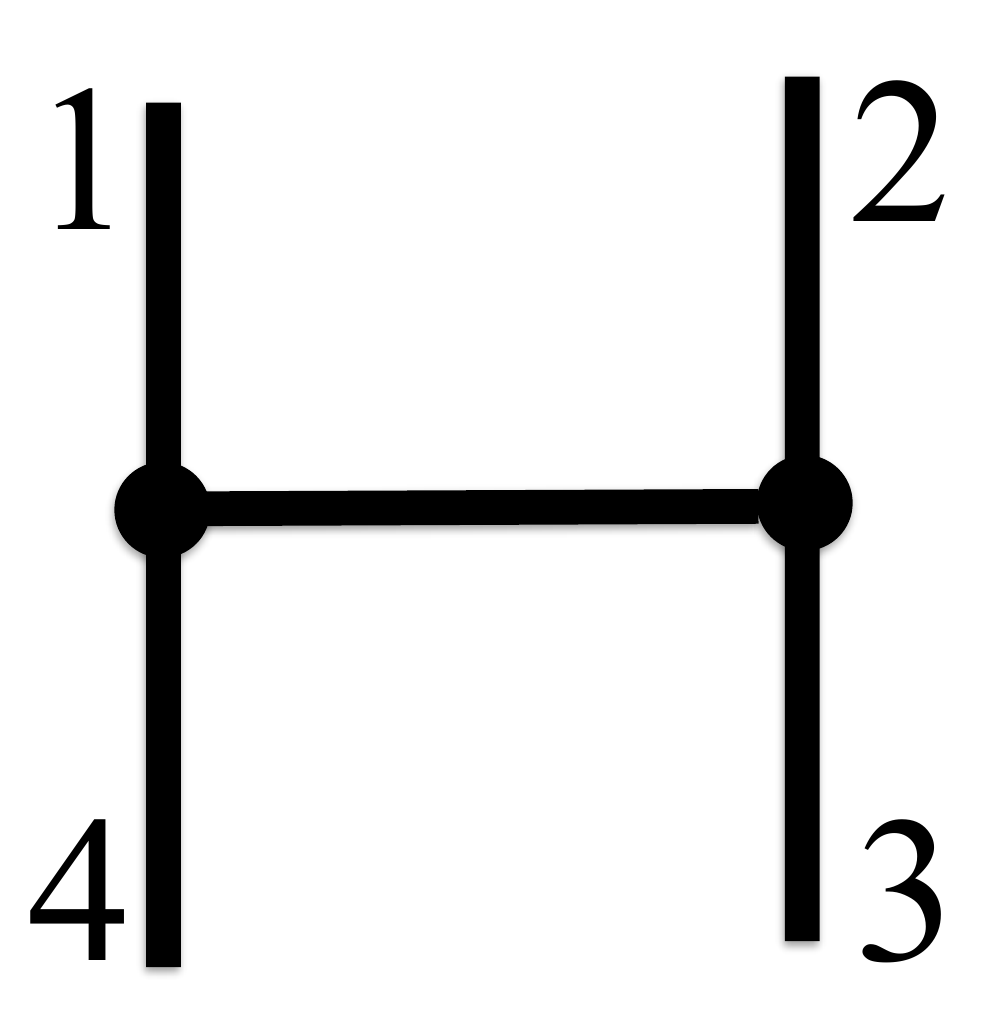}}}~~
+
~~
\raisebox{-.4\height}{\scalebox{0.1}{\includegraphics{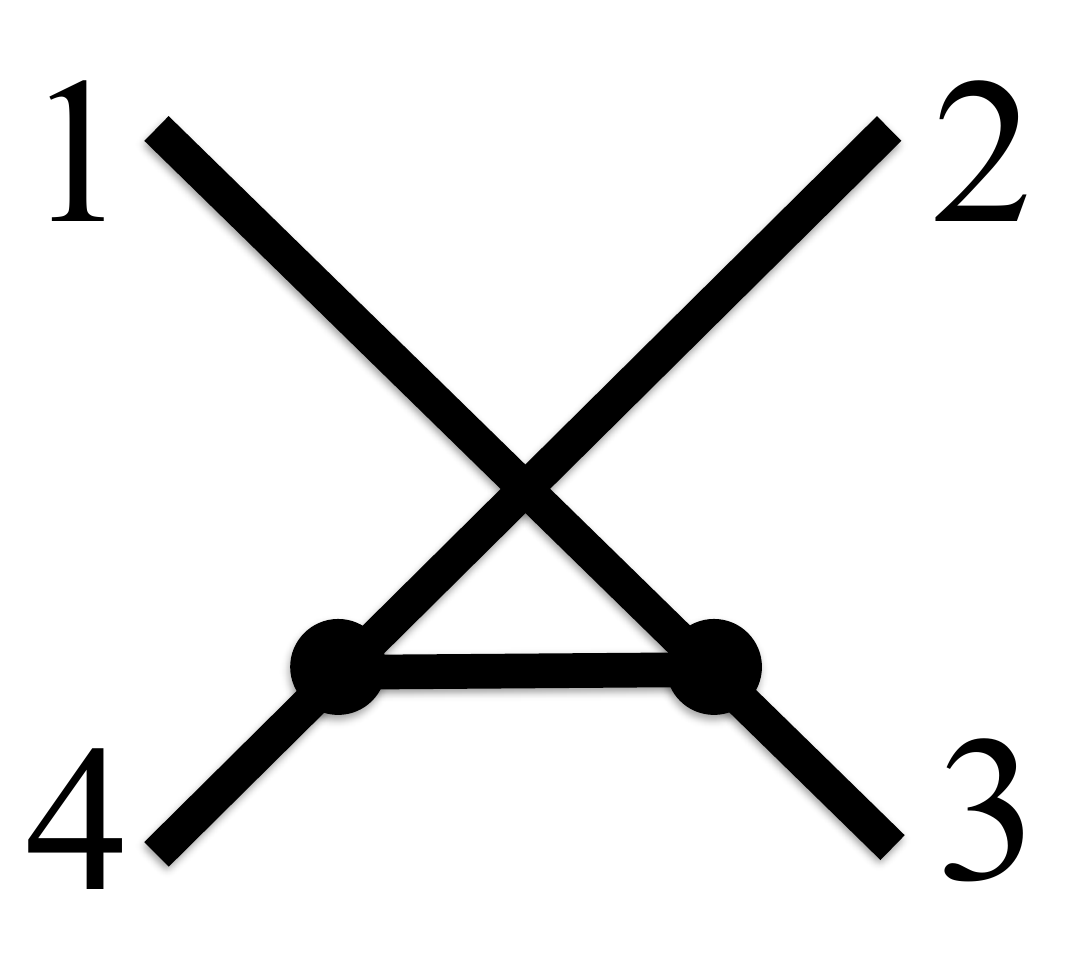}}}~~.
}
As $f$ is a weight system,
$f(\AS\cdot H)=0$ for each $H\in\FF_3$ and $f(\IHX\cdot H)=0$ for each $H\in\FF_4$.
Hence the Claim implies that $\widehat p(\AS)(c)=0$ and $\widehat p(\IHX)(c)=0$.
Therefore, $c=c_{\frak{g}}$ for some metric Lie algebra $\frak{g}$
(cf.\ \rf{12ok14c}).

We show that $\frak{g}$ is reductive.
For this it suffices
to show that the orthogonal complement $Z(\frak{g})^{\perp}$ of the center
$Z(\frak{g})$ of $\frak{g}$ is semisimple (as then $Z(\frak{g})\cap Z(\frak{g})^{\perp}=0$,
so $Z(\frak{g})$ is nondegenerate).

Suppose to the contrary that $Z(\frak{g})^{\perp}$ contains a nonzero abelian ideal $I$.
We can assume that $I$ is a minimal nonzero ideal.
Then $I\subseteq I^{\perp}$, since, by the minimality of $I$, either $[\frak{g},I]=0$, hence
$I\subseteq Z(\frak{g})\subseteq I^{\perp}$,
or $[\frak{g},I]=I$, hence $\langle I,I\rangle=\langle [\frak{g},I],I\rangle=\langle \frak{g},[I,I]\rangle=0$.

So $I+Z(\frak{g})\subseteq I^{\perp}$.
This implies that we can choose a subspace $A$ of $I^{\perp}$ with $I\cap A=0$ and
$I+A=I^{\perp}$ such that
\dyy{18ma14a}{
Z(\frak{g})=(I\cap Z(\frak{g}))+(A\cap Z(\frak{g})).
}
Then $A$ is nondegenerate, since
$A\cap A^{\perp}=A\cap I^{\perp}\cap A^{\perp}=A\cap(I+A)^{\perp}=A\cap I=0$.

So also $A^{\perp}$ is nondegenerate.
As $I\subseteq A^{\perp}$ and $\dim(A^{\perp})=2\dim(I)$,
there exists a self-orthogonal subspace $C$ of $A^{\perp}$ with
$I\cap C=0$ and $I+C=A^{\perp}$.
Then $\dim(C)=\dim(I)$ (as $\dim(A)=n-2\dim(I)$).

Now define, for any nonzero $\alpha\in\oC$, $\varphi_{\alpha}:\frak{g}\to \frak{g}$
by
\dyy{14ma14a}{
\varphi_{\alpha}(x)=
\begin{cases}
\alpha^{-1}x&\text{ if $x\in I$,}\\
\alpha x&\text{ if $x\in C$,}\\
x&\text{ if $x\in A$.}
\end{cases}
}
So $\varphi_{\alpha}\in O(n)$.

Let $\pi_A$ denote the orthogonal projection $\frak{g}\to A$.
Then
\dyy{27ma14a}{
\lim_{\alpha\to 0}c_\frak{g}\cdot\varphi_{\alpha}^{\otimes 3}=c_\frak{g}\cdot\pi_A^{\otimes 3},
}
where $\cdot$ denotes the standard inner product on $(\oC^n)^{\otimes 3}$.
To prove \rf{27ma14a}, choose $x,y,z\in I\cup C\cup A$.
If $x,y,z\in A$, then for each nonzero $\alpha\in\oC$:
\dyyz{
 c_\frak{g}\cdot\varphi_{\alpha}^{\otimes 3}(x\otimes y\otimes z)=
 c_\frak{g}\cdot(x\otimes y\otimes z)=
 c_\frak{g}\cdot\pi_A^{\otimes 3}(x\otimes y\otimes z).
}
If not all of $x,y,z$ belong to $A$, let
$k$ be the number of $x,y,z$ belonging to $I$ minus the number of $x,y,z$
belonging to $C$.
Then
\dyy{27ma14b}{
\lim_{\alpha\to 0} c_\frak{g}\cdot\varphi_{\alpha}^{\otimes 3}(x\otimes y\otimes z)=
\lim_{\alpha\to 0}\alpha^{-k} c_\frak{g}\cdot(x\otimes y\otimes z)=
0.
}
The last equality follows from the fact that if $k\geq 0$, then we may assume
(by symmetry) that
$x\in I$ and $z\not\in C$, so $z\in I+A=I^{\perp}$.
Hence $ c_\frak{g}\cdot(x\otimes y\otimes z)=\langle [x,y],z\rangle=0$, as $[x,y]\in I$.

This proves \rf{27ma14a}.
Hence, as $O(n)\cdot c_\frak{g}$ is closed, there exists $\varphi\in O(n)$
such that
\dyy{23ok14a}{
 c_\frak{g}\cdot\pi_A^{\otimes 3}= c_\frak{g}\cdot\varphi^{\otimes 3}.
}
This implies
\dyy{27ma14d}{
\varphi(I+C+Z(\frak{g}))\subseteq Z(\frak{g}).
}
To see this, by \rf{18ma14a}, $I+C+Z(\frak{g})=I+C+(A\cap Z(\frak{g}))$.
Now choose $x\in I\cup C\cup (A\cap Z(\frak{g}))$.
Then for all $y,z\in \frak{g}$, using \rf{23ok14a}:
\dyy{27ma14c}{
\langle[\phi(x),\phi(y)],\phi(z)\rangle=
\langle[\pi_A(x),\pi_A(y)],\pi_A(z)\rangle=
0.
}
Indeed, if $x\in I+C=A^{\perp}$ then $\pi_A(x)=0$.
If $x\in A\cap Z(\frak{g})$, then $\pi_A(x)=x$, and \rf{27ma14c} follows as $x\in Z(\frak{g})$.
As \rf{27ma14c} holds for all $y,z\in\frak{g}$, $\varphi(x)$ belongs to $Z(\frak{g})$.

So we have \rf{27ma14d}, which implies $\dim(I+C+Z(\frak{g}))\leq\dim(Z(\frak{g}))$,
so $C\subseteq Z(\frak{g})$, hence $C=0$, as $C\cap Z(\frak{g})\subseteq C\cap (I+A)=0$.
Therefore,
$\dim(I)=\dim(C)=0$, contradicting $I\neq 0$.
Concluding, $\frak{g}$ is reductive.
\bx

\sect{8ok14j}{Proof of the theorem: uniqueness of $\frak{g}$}

We first show uniqueness if the metrics are the Killing forms.

\prop{25ma14b}{
Let $\frak{g}$ and $\frak{g}'$ be complex semisimple Lie algebras with their Killing forms as metrics.
If $\varphi_\frak{g}=\varphi_{\frak{g}'}$ then $\frak{g}= \frak{g}'$.
}

\pf
As $\phi_{\frak{g}}=\phi_{\frak{g}'}$, we know
$\dim(\frak{g})=\varphi_\frak{g}(\bgcirc)=\varphi_{\frak{g}'}(\bgcirc)=\dim(\frak{g}')$.
Let $\frak{h}$ and $\frak{h}'$ be real compact forms in $\frak{g}$ and $\frak{g}'$ respectively.
Since the Killing forms are negative definite on $\frak{h}$ and $\frak{h}'$, we can assume that the 
inner product spaces underlying $\frak{h}$ and $\frak{h}'$ both are $\oR^n$ with
standard negative definite inner product,
and that $c_{\frak{g}}$ and $c_{\frak{g}'}$ belong to $((\oR^n)^{\otimes 3})^{C_3}$.

Suppose $\frak{g}\neq\frak{g}'$.
Hence the orbits $O(n,\oR)\cdot c_\frak{g}$ and $O(n,\oR)\cdot c_{\frak{g}'}$ are disjoint
compact subsets of $((\oR^n)^{\otimes 3})^{C_3}$.
By the Stone-Weierstrass theorem, there exists a real-valued polynomial
$q$ on $((\oR^n)^{\otimes 3})^{C_3}$ such that $q(x)\leq 0$ for each $x\in O(n,\oR)\cdot c_\frak{g}$
and $q(x)\geq 1$ for each $x\in O(n,\oR)\cdot c_{\frak{g}'}$.
Applying the Reynolds operator, we may assume that $q$ is $O(n,\oR)$-invariant.
By the first fundamental theorem of invariant theory,
$q$ belongs to the algebra generated by $\{p(G)\mid G$ 3-graph$\}$.
However, $p(G)(c_{\frak{g}})=\varphi_{\frak{g}}(G)=\varphi_{\frak{g}'}(G)=p(G)(c_{\frak{g}'})$
for each $3$-graph $G$.
So $q(c_{\frak{g}})=q(c_{\frak{g}'})$, contradicting $q(c_{\frak{g}})\leq 0$ and $q(c_{\frak{g}'})\geq 1$.
\obx

For each complex metric Lie algebra ${\frak{g}}$ of positive dimension, define
\dyyz{
\varphi'_{\frak{g}}:=\frac{1}{\dim(\frak{g})}\varphi_{\frak{g}}.
}
From Proposition \ref{25ma14b} we derive the next proposition.

\prop{25ma14c}{
Let $\frak{g}$ and $\frak{g}'$ be complex simple metric Lie algebras.
If $\varphi'_\frak{g}=\varphi'_{\frak{g}'}$ then $\frak{g}= \frak{g}'$.
}

\pf
Let $B$ and $B'$ denote the bilinear forms associated with $\frak{g}$ and $\frak{g}'$, respectively,
and let $K$ and $K'$ be the Killing forms of $\frak{g}$ and $\frak{g}'$, respectively.
Since $\frak{g}$ and $\frak{g}'$ are simple, there are nonzero $\alpha,\alpha'\in\oC$
such that $B=\alpha K$ and $B'=\alpha' K'$.
Then
\dyyz{
\varphi_{\frak{g},B}(\Thetagraaf)=\alpha^{-1} \varphi_{\frak{g},K}(\Thetagraaf)=
-\alpha^{-1}K^{\otimes 3}( c_{\frak{g},K}, c_{\frak{g},K})
=
\alpha^{-1}\dim(\frak{g}),
}
and similarly $\varphi_{\frak{g}',B'}(\Thetagraaf)={\alpha'}^{-1}\dim(\frak{g}')$.
Since $\varphi'_{\frak{g},B}(\Thetagraaf)=\varphi'_{\frak{g}',B'}(\Thetagraaf)$, this implies
$\alpha=\alpha'$.
So $\varphi'_{\frak{g},K}=\varphi'_{\frak{g}',K'}$, hence we can assume that $\alpha=1$,
so $B=K$ and $B'=K'$.

Now let $\widetilde{\frak{g}}$ be the direct sum of $\dim(\frak{g}')$ copies of $\frak{g}$.
Similarly, let $\widetilde{\frak{g}'}$ be the direct sum of $\dim(\frak{g})$ copies of $\frak{g}'$.
So $\dim\widetilde{\frak{g}}=\dim\widetilde{\frak{g}}'$, and for each $3$-graph $G$, as
$\phi'_{\frak{g}}=\phi'_{\frak{g}'}$ and as $G$ is connected:
\dyyz{
\varphi_{\widetilde{\frak{g}}}(G)=\dim(\frak{g}')\varphi_\frak{g}(G)=\dim(\frak{g})\varphi_{\frak{g}'}(G)=\varphi_{\widetilde{\frak{g}}'}(G).
}
Hence by Proposition \ref{25ma14b}, $\widetilde{\frak{g}}=\widetilde{\frak{g}}'$, and so
$\frak{g}=\frak{g}'$.
\obx

We note that also if $\frak{g}$ is a complex $1$-dimensional metric Lie algebra and $\frak{g}'$ is a
complex simple metric Lie algebra, then $\varphi'_{\frak{g}}\neq\varphi'_{\frak{g}'}$,
since $\varphi'_{\frak{g}}(\Thetagraaf)=0$ while $\varphi'_{\frak{g}'}(\Thetagraaf)\neq 0$.
This and Proposition \ref{25ma14c} is used to prove the last proposition,
which settles the theorem.

\prop{3ok14b}{
Let $\frak{g}$ and $\frak{g}'$ be complex reductive metric Lie algebras.
If $\varphi_\frak{g}=\varphi_{\frak{g}'}$ then $\frak{g}= \frak{g}'$.
}

\pf
As $\frak{g}$ and $\frak{g}'$ are reductive, we can write
\dyz{
$\dps \frak{g}=\bigoplus_{i=1}^m\frak{g}_i$
and
$\dps \frak{g}'=\bigoplus_{j=1}^{m'}\frak{g}'_j$,
}
where each $\frak{g}_i$ and $\frak{g}'_j$ is either simple or $1$-dimensional.
Then, since $3$-graphs are connected,
\dyy{21ma14b}{
\sum_{i=1}^m\varphi_{\frak{g}_i}
=
\varphi_{\frak{g}}
=
\varphi_{\frak{g}'}
=
\sum_{j=1}^{m'}\varphi_{\frak{g}'_j}.
}
So we can assume that $\frak{g}_i\neq\frak{g}'_j$ for all $i\in[m]$ and $j\in[m']$.
Hence, by Proposition \ref{25ma14c} and the remark thereafter, there exist finitely many 3-graphs
$G_1,\ldots,G_k$ such that for all $i\in[m]$ and $j\in[m']$ there exists $t\in[k]$ with
$\varphi'_{\frak{g}_i}(G_t)\neq\varphi'_{\frak{g}'_j}(G_t)$.
That is, for each $i\in [m]$ and $j\in[m']$, the following
vectors ${\mathbf y}_i$, ${\mathbf z}_j\in\oC^k$:
\dyz{
$\dps{\mathbf y}_i:=(\varphi'_{\frak{g}_i}(G_1),\ldots,\varphi'_{\frak{g}_i}(G_k))$
~~~and~~~
$\dps{\mathbf z}_j:=(\varphi'_{\frak{g}'_j}(G_1),\ldots,\varphi'_{\frak{g}'_j}(G_k))$
}
are distinct.
So there exists a polynomial $q\in\oC[x_1,\ldots,x_k]$ such that
$q({\mathbf y}_i)=0$ for each $i=1,\ldots,m$ and
$q({\mathbf z}_j)=1$ for each $j=1,\ldots,m'$.
Now set $\gamma:=q(G_1,\ldots,G_k)$, taking formal linear sums of 3-graphs and
applying the following composition 
of $3$-graphs $G$ and $H$ as product ([8]):
take the disjoint union of $G$ and $H$, choose an edge
$uv$ of $G$ and an edge $u'v'$ of $H$, and replace them by $uu'$ and $vv'$.
Let $F$ be the 3-graph thus arising.
Then for any complex simple or 1-dimensional metric Lie algebra
$\frak{g}$: $\varphi'_\frak{g}(F)=\varphi'_\frak{g}(G)\varphi'_\frak{g}(H)$,
independently of the choice of $uv$ and $u'v'$
(see Proposition 7.18 in [6]).

We extend each $\varphi'_{\frak{g}_i}$ and $\varphi'_{\frak{g}'_j}$ linearly to $\gamma$.
Then $\varphi'_{\frak{g}_i}(\gamma)=q({\mathbf y}_i)=0$ for each $i=1,\ldots,m$
while $\varphi'_{\frak{g}'_j}(\gamma)=q({\mathbf z}_j)=1$ for each $j=1,\ldots,m'$.
Hence $\varphi_{\frak{g}_i}(\gamma)=0$ for each $i=1,\ldots,m$
and $\varphi_{\frak{g}'_j}(\gamma)=\dim(\frak{g}'_j)$ for each $j=1,\ldots,m'$.
Therefore, by \rf{21ma14b}, $m'=0$.
Similarly, $m=0$.
\bx

\section*{References}\label{REF}
{\small
\begin{itemize}{}{
\setlength{\labelwidth}{4mm}
\setlength{\parsep}{0mm}
\setlength{\itemsep}{1mm}
\setlength{\leftmargin}{5mm}
\setlength{\labelsep}{1mm}
}
\item[\mbox{\rm[1]}] S. Axelrod, I.M. Singer, 
Chern-Simons perturbation theory {II},
{\em Journal of Differential Geometry} 39 (1994), 173--213.

\item[\mbox{\rm[2]}] D. Bar-Natan, 
{\em Perturbative Aspects of the Chern-Simons Topological Quantum Field Theory},
Ph.D. Thesis, Harvard University, Cambridge, Mass., 1991.

\item[\mbox{\rm[3]}] D. Bar-Natan, 
On the Vassiliev knot invariants,
{\em Topology} 34 (1995) 423--472.

\item[\mbox{\rm[4]}] D. Bar-Natan, 
Lie algebras and the four color theorem,
{\em Combinatorica} 17 (1997) 43--52.

\item[\mbox{\rm[5]}] M. Brion, 
Introduction to actions of algebraic groups,
{\em Les cours du C.I.R.M.} 1 (2010) 1--22.

\item[\mbox{\rm[6]}] S. Chmutov, S. Duzhin, J. Mostovoy, 
{\em Introduction to Vassiliev Knot Invariants},
Cambridge University Press, Cambridge, 2012.

\item[\mbox{\rm[7]}] J. Draisma, D. Gijswijt, L. Lov\'asz, G. Regts, A. Schrijver, 
Characterizing partition functions of the vertex model,
{\em Journal of Algebra} 350 (2012) 197--206.

\item[\mbox{\rm[8]}] S.V. Duzhin, A.I. Kaishev, S.V. Chmutov, 
The algebra of $3$-graphs,
{\em Proceedings of the Steklov Institute of Mathematics} 221 (1998) 157--186.

\item[\mbox{\rm[9]}] M.H. Freedman, L. Lov\'asz, A. Schrijver, 
Reflection positivity, rank connectivity, and homomorphisms of graphs,
{\em Journal of the American Mathematical Society} 20 (2007) 37--51.

\item[\mbox{\rm[10]}] R. Goodman, N.R. Wallach, 
{\em Symmetry, Representations, and Invariants},
Springer, Dordrecht, 2009.

\item[\mbox{\rm[11]}] P. de la Harpe, V.F.R. Jones, 
Graph invariants related to statistical mechanical models:
examples and problems,
{\em Journal of Combinatorial Theory, Series B} 57 (1993) 207--227.

\item[\mbox{\rm[12]}] M. Kontsevich, 
Feynman diagrams and low-dimensional topology,
in: {\em First European Congress of Mathematics Volume II},
Birkha\"user, Basel, 1994, pp. 97--121.

\item[\mbox{\rm[13]}] T. Murphy, 
On the tensor system of a semisimple Lie algebra,
{\em Proceedings of the Cambridge Philosophical Society
(Mathematical and Physical Sciences)}
71 (1972) 211--226.

\item[\mbox{\rm[14]}] R. Penrose, 
Applications of negative dimensional tensors,
in: {\em Combinatorial Mathematics and Its Applications}
(D.J.A. Welsh, ed.),
Academic Press, London, 1971, pp. 221--244.

\item[\mbox{\rm[15]}] B. Szegedy, 
Edge coloring models and reflection positivity,
{\em Journal of the American Mathematical Society}
20 (2007) 969--988.

\item[\mbox{\rm[16]}] P. Vogel, 
Algebraic structures on modules of diagrams,
{\em Journal of Pure and Applied Algebra} 215 (2011) 1292--1339.

\item[\mbox{\rm[17]}] H. Weyl, 
{\em The Classical Groups --- Their Invariants and Representations},
Princeton University Press, Princeton, New Jersey, 1946.

\end{itemize}
}

\end{document}